\documentclass[11pt]{article}
\usepackage{latexsym}

\newcommand{\eproof}{\mbox{\ }\hfill $\Box$ \par \vskip 10pt}

\newtheorem{Theorem}{Theorem}[section]
\newtheorem{lemma}[Theorem]{Lemma}
\newtheorem{prop}[Theorem]{Proposition}

\begin{document}

\title{Resolvent estimates for the magnetic Schr\"odinger operator}

\author{{\sc Georgi Vodev}}

\date{}

\maketitle

\noindent
{\bf Abstract.} We prove optimal high-frequency resolvent estimates for self-adjoint operators of the form
$G=-\Delta+ib(x)\cdot\nabla+i\nabla\cdot b(x)+V(x)$ on $L^2({\bf R}^n)$, $n\ge 3$, where $b(x)$ and $V(x)$ are large magnetic and electric potentials, respectively.

\setcounter{section}{0}
\section{Introduction and statement of results}

The purpose of this work is to extend the recent results in \cite{kn:CCV1}, \cite{kn:CCV2} to 
a larger class (most probably optimal) of magnetic potentials.
In the present paper we study the high-frequency behaviour of the resolvent of self-adjoint operators of the form
$$G=-\Delta+ib(x)\cdot\nabla+i\nabla\cdot b(x)+V(x)\quad{\rm on}\quad L^2({\bf R}^n),\,\, n\ge 3,$$
where $b=(b_1,...,b_n)\in L^\infty({\bf R}^n;{\bf R}^n)$ is a magnetic potential
and $V\in L^\infty({\bf R}^n;{\bf R})$ is an electric potential. Introduce the polar coordinates $r=|x|$, $w=\frac{x}{|x|}\in{\bf S}^{n-1}$. 
We suppose that $b(x)=b^L(x)+b^S(x)$, $V(x)=V^L(x)+V^S(x)$ with long-range parts $b^L$ and $V^L$ belonging to $C^1([r_0,+\infty))$, 
$r_0\gg 1$, with respect to the radial variable $r$ and satisfying the conditions
$$\left|V^L(rw)\right|\le C,\eqno{(1.1)}$$
$$\partial_r V^L(rw)\le Cr^{-1-\delta},\eqno{(1.2)}$$
$$\left|\partial_r^kb^L(rw)\right|\le Cr^{-k-\delta},\quad k=0,1,\eqno{(1.3)}$$
for all $r\ge r_0$, $w\in {\bf S}^{n-1}$, with some constants $C,\delta>0$. The short-range parts satisfy
$$\left|b^S(x)\right|+\left|V^S(x)\right|\le C\langle x\rangle^{-1-\delta}.\eqno{(1.4)}$$ 
Our main result is the following 

\begin{Theorem}  Under the conditions (1.1)-(1.4), for every
$s>1/2$ there exist constants $C,\lambda_0>0$ so that for $\lambda\ge\lambda_0$, $0<\varepsilon\le 1$, 
$|\alpha_1|,|\alpha_2|\le 1$, we have the estimate
 $$\left\|\langle x\rangle^{-s}\partial_{x}^{\alpha_1}
\left(G-\lambda^2\pm i\varepsilon\right)^{-1}\partial_{x}^{\alpha_2}\langle x\rangle^{-s}\right\|_{L^2\to L^2}\le
C\lambda^{|\alpha_1|+|\alpha_2|-1}.\eqno{(1.5)}$$
\end{Theorem} 

Note that this kind of resolvent estimates play an important role in proving uniform local energy decay, 
 dispersive, smoothing and Strichartz estimates for solutions to the corresponding wave and
Schr\"odinger equations (e.g. see \cite{kn:CCV1}, \cite{kn:CCV3}, \cite{kn:EGS}).

Theorem 1.1 is proved in \cite{kn:CCV1} assuming additionally that $b^S(x)$ is continuous with respect to the radial variable $r$ uniformly in $w$.
In the case $b^L\equiv 0$, $V^L\equiv 0$ and $b^S$, $V^S$ satisfying (1.4), the estimate (1.5)
is proved in \cite{kn:EGS} under the extra assumption that $b(x)$ is continuous in $x$. 
In fact, no continuity of the magnetic potential is needed in order to have (1.5) as shown in \cite{kn:CCV2}. 
Instead, it was supposed in \cite{kn:CCV2} that ${\rm div}\,b^L$ and ${\rm div}\,b^S$ exist as functions in $L^\infty$. This assumption allows to
conclude that the perturbation (which is a first-order differential operator) sends
the Sobolev space $H^1$ into $L^2$, a fact used in an essential way in \cite{kn:CCV2}. 
Note also the work \cite{kn:KT} where it is shown that 
in the case $b^L\equiv 0$, $V^L\equiv 0$ and $b^S$, $V^S$ satisfying (1.4) the operator $G$ 
 has no strictly positive eigenvalues. It follows from Theorem 1.1 that in the more general case when
the long-range parts are not identically zero the operator $G$ 
 has no strictly positive eigenvalues outside some compact interval.
 
There are two main difficulties in proving the above theorem. The first one is that under our assumptions the commutator of the gradient
and the magnetic potential is not an $L^\infty$ function. Consiquently, the perturbation does not send
the Sobolev space $H^1$ into $L^2$. Instead, it is bounded from $H^1$ into $H^{-1}$. Secondly, the magnetic potential is large
and therefore it is hard to apply perturbation arguments similar to those used in \cite{kn:CCV1}. Thus, 
 to prove Theorem 1.1 we first observe that (1.5) is equivalent to a semi-classical a priori estimate on weighted 
 Sobolev spaces (see (2.10) below). Furthermore, we derive this a priori estimate from a semi-classical Carleman estimate
 on weighted Sobolev spaces (see (2.7) below) with a suitably chosen phase function independent of the semi-classical parameter.
 To get this Carleman estimate we first prove a semi-classical Carleman estimate
 on weighted Sobolev spaces for the long-range part of the operator (see Theorem 2.1 below) and we then apply a perturbation argument.
 Note that the estimate (2.1) is valid for any phase function $\varphi(r)\in C^2({\bf R})$ such that its first derivative $\varphi'(r)$ is of compact support and non-negative. The main feature of our Carleman estimate is that it is uniform with respect to the phase function
 $\varphi$ (that is, the constant $C_1$ does not depend on $\varphi$) and the weight in the right-hand side is smaller than the usual one
 (that is, $\left(\langle x\rangle^{-2s}+\varphi'(|x|)\right)^{-1/2}$ instead of $\langle x\rangle^{s}$). Thus we can make this weight
 small on an arbitrary compact set by choosing the phase function properly. 
 Moreover, in the right-hand side
 we have the better semi-classical Sobolev $H^{-1}$ norm instead of the $L^2$ one, which is crucial for the application we make here. 
 Note also that Carleman estimates similar to (2.1) and (2.7) have been recently proved in \cite{kn:D} for operators of the form
 $-h^2\Delta+V(x,h)$, where $V$ is a real-valued long-range potential which is $C^1$ with respect to the radial variable $r$.
 There are, however, several important differences between the Carleman estimates in \cite{kn:D} and ours. First, the phase function in \cite{kn:D}
 is of the form $\varphi=\varphi_1(r)/h$, where $\varphi_1$ does not depend on $h$ and must satisfy some conditions.
 Thus, the Carleman estimates in \cite{kn:D} lead to the conclusion that the resolvent in that case is bounded by $e^{C/h}$, $C>O$ being a constant. 
 Secondly, in \cite{kn:D} the Carleman estimates are not uniform with respect to the phase function and the norm in the right-hand side is $L^2$ (and not $H^{-1}$). Finally, the operator in \cite{kn:D} does not contain a magnetic potential.
 
 To prove Theorem 2.1 we make use of methods originating from \cite{kn:CV}. Note that in \cite{kn:CV} the high-frequency behavior
of the resolvent of operators of the form $-\Delta_g+V$ is studied, where $V$ is a real-valued scalar potential and $\Delta_g$ is the negative
Laplace-Beltrami operator on unbounded Riemannian manifolds as for example asymptotically Euclidean and hyperbolic ones. Similar techniques
have been also used in \cite{kn:RT} where actually all ranges of frequencies are covered. In these two papers, however, no perturbations
by magnetic potentials are studied.

\section{Proof of Theorem 1.1}

Set $h=\lambda^{-1}$, $P(h)=h^2G$, $\widetilde b(x,h)=hb(x)$, $\widetilde b^L(x,h)=h\chi(|x|)b^L(x)$, $\widetilde b^S(x,h)=
\widetilde b(x,h)-\widetilde b^L(x,h)$, 
$\widetilde V(x,h)=h^2V(x)$, $\widetilde V^L(x,h)=h^2\chi(|x|)V^L(x)$, $\widetilde V^S(x,h)=
\widetilde V(x,h)-\widetilde V^L(x,h)$, where 
$\chi\in C^\infty({\bf R})$, 
$\chi(r)=0$ for $r\le r_0+1$, $\chi(r)=1$ for $r\ge r_0+2$. 
Throughout this paper $H^1({\bf R}^n)$ will denote the Sobolev space equipped with the semi-classical norm
$$\left\|u\right\|_{H^1}^2=\sum_{0\le|\alpha|\le 1}\left\|{\cal D}_x^\alpha u\right\|_{L^2}^2$$
where ${\cal D}_x=ih\partial_x$. Furthermore, $H^{-1}$ will denote the dual space of $H^1$ with respect to
the scalar product $\langle\cdot,\cdot\rangle_{L^2}$ with a norm
$$\left\|v\right\|_{H^{-1}}=\sup_{0\neq u\in H^1}\frac{\left|\langle u,v\rangle_{L^2}\right|}{\left\|u\right\|_{H^1}}.$$
Let $\rho\in C^{\infty}({\bf R})$ be a function independent of $h$ such that  $0\le\rho\le 1$ and 
 $\rho(\sigma)=1$ for $\sigma\le 0$, $\rho(\sigma)=0$ for $\sigma\ge 1$. 
Define the function $\varphi(r)\in C^{\infty}({\bf R})$ as follows: $\varphi(0)=0$ and
$$\varphi'(r)=\tau\rho(r-A)$$
where $\tau, A\ge 1$ are parameters independent of $h$ to be fixed later on. 
Introduce the operator
$$P^L(h)=-h^2\Delta+ih\widetilde b^L(x,h)\cdot\nabla+ih\nabla\cdot \widetilde b^L(x,h)+\widetilde V^L(x,h)$$
and set
$$P^L_\varphi(h)=e^\varphi P^L(h)e^{-\varphi},$$
$$P_\varphi(h)=e^\varphi P(h)e^{-\varphi}$$ $$=P^L_\varphi(h)+ih\widetilde b^S(x,h)\cdot\nabla+ih\nabla\cdot\widetilde b^S(x,h) 
-2ih\widetilde b^S(x,h)\cdot\nabla\varphi+\widetilde V^S(x,h),$$

$$\mu(x)=\sqrt{\langle x\rangle^{-2s}+\varphi'(|x|)}.$$

In what follows in this section we will show that Theorem 1.1 follows from the following

\begin{Theorem} Suppose (1.1), (1.2), (1.3) fulfilled and let $\frac{1}{2}<s<\frac{1+\delta}{2}$. Then, 
 we have the a priori estimate
$$\left\|\langle x\rangle^{-s} f\right\|_{H^1}\le\frac{C_1}{h}\left\|\mu^{-1}\left(P_\varphi^L(h)-1\pm i\varepsilon\right)f\right\|_{H^{-1}}
+C_2\left(\frac{\varepsilon}{h}\right)^{1/2}\left\|f\right\|_{L^2}\eqno{(2.1)}$$
for $0<\varepsilon\le 1$, $0<h\le h_0(\tau,A)\ll 1$, with a constant $C_1>0$ independent of $f$, 
$\varepsilon$, $h$, $\tau$, $A$, and a constant $C_2>0$ independent of 
$f$, $\varepsilon$, $h$.    
\end{Theorem}

Let us first see that (2.1) implies the estimate
$$\left\|\langle x\rangle^{-s}f\right\|_{H^1}\le\frac{2C_1}{h}\left\|\langle x\rangle^{s}
\left(P_\varphi(h)-1\pm i\varepsilon\right)f\right\|_{H^{-1}}
+2C_2\left(\frac{\varepsilon}{h}\right)^{1/2}\left\|f\right\|_{L^2}.\eqno{(2.2)}$$
Using that $\mu(x)\ge\tau^{1/2}$ for $|x|\le A$ and $\mu(x)\ge\langle x\rangle^{-s}$ for $|x|\ge A+1$ together with
the condition (1.4), we get (for $0<s-\frac{1}{2}\ll 1$)
$$\langle x\rangle^s\mu(x)^{-1}\left(\left|\widetilde b^S(x,h)\right|+\left|\widetilde V^S(x,h)\right|\right)\le
Ch\left(\tau^{-1/2}+A^{2s-1-\delta}\right),\eqno{(2.3)}$$
$$\langle x\rangle^s\mu(x)^{-1}\left|\widetilde b^S(x,h)\right|\left|\nabla\varphi\right|\le
O_{\tau,A}(h).\eqno{(2.4)}$$
By (2.3) and (2.4),
$$\left\|\mu^{-1}\left(P_\varphi(h)-P_\varphi^L(h)\right)\langle x\rangle^s\right\|_{H^1\to H^{-1}}\le
Ch\left(\tau^{-1/2}+A^{2s-1-\delta}+O(h)\right).\eqno{(2.5)}$$
By (2.1) and (2.5),
$$\left\|\langle x\rangle^{-s} f\right\|_{H^1}\le\frac{C_1}{h}\left\|\mu^{-1}\left(P_\varphi(h)-
1\pm i\varepsilon\right)f\right\|_{H^{-1}}$$ $$+\frac{C_1}{h}\left\|\mu^{-1}
\left(P_\varphi(h)-P_\varphi^L(h)\right)f\right\|_{H^{-1}}
+C_2\left(\frac{\varepsilon}{h}\right)^{1/2}\left\|f\right\|_{L^2}$$ $$
\le\frac{C_1}{h}\left\|\langle x\rangle^{s}\left(P_\varphi(h)-
1\pm i\varepsilon\right)f\right\|_{H^{-1}}$$ $$+C\left(\tau^{-1/2}+A^{2s-1-\delta}+O(h)\right)
\left\|\langle x\rangle^{-s}f\right\|_{H^{1}}
+C_2\left(\frac{\varepsilon}{h}\right)^{1/2}\left\|f\right\|_{L^2}.\eqno{(2.6)}$$
Taking now $\tau^{-1}$, $A^{-1}$ and $h$ small enough, we can absorb the second term in the right-hand side of (2.6)
to obtain (2.2).

Applying (2.2) with $f=e^\varphi g$ we obtain the following Carleman estimate
$$\left\|\langle x\rangle^{-s}e^\varphi g\right\|_{H^1}\le\frac{2C_1}{h}\left\|\langle x\rangle^{s}e^\varphi 
\left(P(h)-1\pm i\varepsilon\right)g\right\|_{H^{-1}}
+2C_2\left(\frac{\varepsilon}{h}\right)^{1/2}\left\|e^\varphi g\right\|_{L^2}.\eqno{(2.7)}$$
Since the function $\varphi$ does not depend on $h$, we deduce from (2.7) the a priori estimate
$$\left\|\langle x\rangle^{-s}g\right\|_{H^1}\le\frac{\widetilde C_1}{h}\left\|\langle x\rangle^{s}
\left(P(h)-1\pm i\varepsilon\right)g\right\|_{H^{-1}}
+\widetilde C_2\left(\frac{\varepsilon}{h}\right)^{1/2}\left\|g\right\|_{L^2}\eqno{(2.8)}$$
with constants $\widetilde C_1,\widetilde C_2>0$ independent of $h$, $\varepsilon$ and $g$.
On the other hand, since the operator $P(h)$ is symmetric on $L^2({\bf R}^n)$, we have
$$\varepsilon\|g\|_{L^2}^2=\mp{\rm Im}\,\left\langle(P(h)-1\pm i\varepsilon)g,g\right\rangle_{L^2}$$ $$
\le \gamma^{-1}h^{-1}\left\|\langle x\rangle^{s}\left(P(h)-1\pm i\varepsilon\right)g\right\|_{H^{-1}}^2
+\gamma h\left\|\langle x\rangle^{-s}g\right\|_{H^1}^2\eqno{(2.9)}$$
for every $\gamma>0$. 
Taking $\gamma$ small enough, independent of $h$, we deduce from (2.8) and (2.9) the  
a priori estimate
$$\left\|\langle x\rangle^{-s}g\right\|_{H^1}\le\frac{C}{h}\left\|\langle x\rangle^{s}
\left(P(h)-1\pm i\varepsilon\right)g\right\|_{H^{-1}}\eqno{(2.10)}$$
with a constant $C>0$ independent of $h$, $\varepsilon$ and $g$.
It is easy to see now that (2.10) implies the resolvent estimate (1.5) for $0<s-\frac{1}{2}\ll 1$ (and hence
for all $s>\frac{1}{2}$).

\section{Proof of Theorem 2.1} 

We will first prove the following

\begin{prop} Under the conditions of Theorem 2.1 we have the estimate
 $$\left\|\langle x\rangle^{-s} f\right\|_{H^1}\le\frac{C_1}{h}
\left\|\mu^{-1}\left(P_\varphi^L(h)-1\pm i\varepsilon\right)f\right\|_{L^2}
+C_2\left(\frac{\varepsilon}{h}\right)^{1/2}\left\|f\right\|_{H^1}\eqno{(3.1)}$$
for every $0<\varepsilon\le 1$, $0<h\le h_0(\tau,A)\ll 1$, with a constant $C_1>0$ independent of $f$, 
$\varepsilon$, $h$, $\tau$, $A$, and a constant $C_2>0$ independent of 
$f$, $\varepsilon$, $h$.  
\end{prop}

{\it Proof.} We pass to the polar coordinates $(r,w)\in {\bf R}^+\times{\bf S}^{n-1}$, $r=|x|$, $w=\frac{x}{|x|}$, and recall that
$L^2({\bf R}^n)\cong L^2({\bf R}^+\times{\bf S}^{n-1},r^{n-1}drdw)$. Denote by $X$ the Hilbert space $L^2({\bf R}^+\times{\bf S}^{n-1},drdw)$.
We also denote by $\|\cdot\|$ and $\langle\cdot,\cdot\rangle$ the norm and the scalar product on $L^2({\bf S}^{n-1})$. We will make use 
of the identity
$$r^{(n-1)/2}\Delta r^{-(n-1)/2}=\partial_r^2+\frac{\widetilde\Delta_w}{r^2},\eqno{(3.2)}$$
where $\widetilde\Delta_w=\Delta_w-\frac{(n-1)(n-3)}{4}$ and $\Delta_w$ denotes the negative Laplace-Beltrami operator on ${\bf S}^{n-1}$. 
 Observe also that
$$r^{\frac{n-1}{2}}\partial_{x_j}r^{-\frac{n-1}{2}}=w_j\partial_r+r^{-1}q_j(w,\partial_w),\eqno{(3.3)}$$
where $w_j=\frac{x_j}{|x|}$ and $q_j$ is a first order differential operator on ${\bf S}^{n-1}$, independent of $r$, anti-symmetric on 
$L^2({\bf S}^{n-1})$. 
It is easy to see that the operators $Q_j(w,{\cal D}_w)=ihq_j(w,\partial_w)$ and $\Lambda_w=-h^2\widetilde\Delta_w\ge 0$ satisfy the estimate
$$\left\|Q_j(w,{\cal D}_w)v\right\|\le C\left\|\Lambda_w^{1/2}v\right\|+Ch\|v\|,\quad \forall v\in H^1({\bf S}^{n-1}),\eqno{(3.4)}$$
with a constant $C>0$ independent of $h$ and $v$.
Set $u=r^{\frac{n-1}{2}}f$,
$${\cal P}^\pm(h)=r^{\frac{n-1}{2}}\left(P_\varphi^L(h)-1\pm i\varepsilon\right)r^{-\frac{n-1}{2}}.$$
Using (3.2) and (3.3) one can easily check that the operator ${\cal P}^\pm(h)$ can be written in the coordinates $(r,w)$ as follows:
 $${\cal P}^\pm(h)={\cal D}_r^2+
\frac{\Lambda_w}{r^2}-1\pm i\varepsilon+\widetilde V^L+W-2ih\varphi'{\cal D}_r$$ 
$$+\sum_{j=1}^n w_j\left(\widetilde b^L_j(rw,h){\cal D}_r+{\cal D}_r\widetilde b^L_j(rw,h)\right)$$ 
 $$+r^{-1}\sum_{j=1}^n\left(\widetilde b^L_j(rw,h)Q_j(w,{\cal D}_w)+Q_j(w,{\cal D}_w)\widetilde b^L_j(rw,h)\right),$$ 
 where we have put ${\cal D}_r=ih\partial_r$ and 
 $$W=-h^2\varphi'(r)^2-h^2\varphi''(r)-2ih\varphi'\sum_{j=1}^nw_j\widetilde b_j^L.$$
 Set
$$\Phi_s(r)=\left\|\langle r\rangle^{-s}u(r,\cdot)\right\|^2+\left\|\langle r\rangle^{-s}{\cal D}_r u(r,\cdot)
\right\|^2+\left\|\langle r\rangle^{-s}r^{-1}\Lambda_w^{1/2}u(r,\cdot)\right\|^2,$$
$$\Psi_s=\left\|\langle r\rangle^{-s}u\right\|^2_{L^2(X)}+\left\|\langle r\rangle^{-s}{\cal D}_r u\right\|^2_{L^2(X)}
+\left\|\langle r\rangle^{-s}r^{-1}\Lambda_w^{1/2}u\right\|^2_{L^2(X)}=\int_0^\infty\Phi_s(r)dr,$$
$$M^\pm(r)=\left\|{\cal P}^\pm(h)u(r,\cdot)\right\|^2,$$
$${\cal M}^\pm=\int_0^\infty \mu^{-2}M^\pm(r)dr,$$
$$N(r)=\left\|u(r,\cdot)\right\|^2+\left\|{\cal D}_ru(r,\cdot)\right\|^2,$$
$${\cal N}=\int_0^\infty N(r)dr,$$
$$E(r)=-\left\langle\left(r^{-2}\Lambda_w-1+\widetilde V^L\right)u(r,\cdot),u(r,\cdot)\right\rangle
+\left\|{\cal D}_ru(r,\cdot)\right\|^2$$
$$-2r^{-1}\sum_{j=1}^n{\rm Re}\,\left\langle\widetilde b^L_j(rw,h)Q_j(w,{\cal D}_w)u(r,\cdot),u(r,\cdot)\right\rangle.$$
In view of (1.1), (1.3) and (3.4), we have
$$E(r)\ge -\left\|r^{-1}\Lambda_w^{1/2}u(r,\cdot)\right\|^2+\frac{1}{2}\left\|u(r,\cdot)\right\|^2+\left\|{\cal D}_ru(r,\cdot)\right\|^2
-O(h)\Phi_{\frac{1+\delta}{2}}(r),\eqno{(3.5)}$$
provided $h$ is taken small enough. Furthermore, 
using that ${\rm Im}\,\left\langle\widetilde b_j^L{\cal D}_ru,{\cal D}_ru\right\rangle=0$ and $Q_j^*=Q_j$, it is easy to check that $E(r)$ satisfies the identity
$$E'(r):=\frac{dE(r)}{dr}=\frac{2}{r}\left\langle r^{-2}\Lambda_wu(r,\cdot),u(r,\cdot)\right\rangle
-\left\langle\frac{\partial\widetilde V^L}{\partial r}u(r,\cdot),u(r,\cdot)\right\rangle$$
 $$-2\sum_{j=1}^n{\rm Re}\,\left\langle \frac{\partial(\widetilde b^L_j(rw,h)/r)}{\partial r}Q_j(w,{\cal D}_w)u(r,\cdot),u(r,\cdot)\right\rangle$$
 $$-2\sum_{j=1}^n{\rm Re}\,\left\langle w_j\frac{\partial\widetilde b^L_j(rw,h)}{\partial r}u(r,\cdot),{\cal D}_ru(r,\cdot)\right\rangle$$ 
  $$+2h^{-1}{\rm Im}\,\left\langle{\cal P}^\pm(h)u(r,\cdot),{\cal D}_ru(r,\cdot)\right\rangle
 \mp 2\varepsilon h^{-1}{\rm Re}\,\left\langle u(r,\cdot),{\cal D}_ru(r,\cdot)\right\rangle$$ 
$$+4\left\langle \varphi'{\cal D}_ru(r,\cdot),{\cal D}_ru(r,\cdot)\right\rangle-
 2h^{-1}{\rm Im}\,\left\langle Wu(r,\cdot),{\cal D}_ru(r,\cdot)\right\rangle.\eqno{(3.6)}$$
    In view of (1.2), (1.3), (3.4) and (3.6), we obtain the inequality
    $$E'(r)\ge \frac{2}{r}\left\|r^{-1}\Lambda_w^{1/2}u(r,\cdot)\right\|^2
+4\varphi'\left\|{\cal D}_ru(r,\cdot)\right\|^2$$ $$
-2h^{-1}\left\|{\cal P}^\pm(h)u(r,\cdot)\right\|\|{\cal D}_r(r,\cdot)\|
-O(h)\Phi_{\frac{1+\delta}{2}}(r)-O(\varepsilon h^{-1})N(r).\eqno{(3.7)}$$
Since $\Phi_{\frac{1+\delta}{2}}(r)\le \Phi_s(r)$ for $\frac{1}{2}<s\le
\frac{1+\delta}{2}$, we obtain from (3.7)
$$E'(r)\ge \frac{2}{r}\left\|r^{-1}\Lambda_w^{1/2}u(r,\cdot)\right\|^2
+4\varphi'\left\|{\cal D}_ru(r,\cdot)\right\|^2$$ $$
-\gamma^{-1}h^{-2}\mu^{-2}M^\pm(r)-\gamma\mu^2\|{\cal D}_r(r,\cdot)\|^2
-O(h)\Phi_s(r)-O(\varepsilon h^{-1})N(r)$$
 $$\ge \frac{2}{r}\left\|r^{-1}\Lambda_w^{1/2}u(r,\cdot)\right\|^2
-\gamma^{-1}h^{-2}\mu^{-2}M^\pm(r)
-O(h+\gamma)\Phi_s(r)-O(\varepsilon h^{-1})N(r)\eqno{(3.8)}$$
for every $0<\gamma\ll 1$. By (3.5) and (3.8),
$$\langle r\rangle^{-2s}\left(E(r)+rE'(r)\right)$$ $$\ge \Phi_s(r)-\gamma^{-1}h^{-2}\mu^{-2}M^\pm(r)
-O(h+\gamma)\Phi_s(r)-O(\varepsilon h^{-1})N(r).\eqno{(3.9)}$$
Integrating (3.8) from $t>0$ to $+\infty$ we get
    $$E(t)=-\int_t^\infty E'(r)dr\le 
O(\gamma^{-1}h^{-2}){\cal M}^\pm+O(\varepsilon h^{-1}){\cal N}+O(h+\gamma)\Psi_s.\eqno{(3.10)}$$
Let $\psi>0$ be a function independent of $h$ and such that $\int_0^\infty\psi(r)dr<\infty$. 
 Multiplying (3.10) by $\psi(t)$ and integrating from $0$ to $+\infty$, we get
    $$\int_0^\infty \psi(r)E(r)dr
\le O(\gamma^{-1}h^{-2}){\cal M}^\pm+O(\varepsilon h^{-1}){\cal N}+O(h+\gamma)\Psi_s.\eqno{(3.11)}$$
Observe now that we have the identity
    $$\int_0^\infty \langle r\rangle^{-2s}\left(E(r)+rE'(r)\right)dr=\int_0^\infty\psi(r)E(r)dr\eqno{(3.12)}$$
    where $\psi(r)=2sr\langle r\rangle^{-2s-1}$. 
 Combining (3.9), (3.11) and (3.12) and taking $\gamma$ and $h$ small enough, we conclude
$$\Psi_s\le O(h^{-2}){\cal M}^\pm+O(\varepsilon h^{-1}){\cal N}.\eqno{(3.13)}$$ 
Clearly, (3.13) implies (3.1).
\eproof

In what follows in this section we will show that (2.1) follows from (3.1) and the following

\begin{lemma} Let $\ell\in{\bf R}$. Then we have the estimate
$$\left\|\mu^{-\ell}\left(P_\varphi^L(h)-i\right)^{-1}\mu^\ell\right\|_{H^{-1}\to H^1}\le C\eqno{(3.14)}$$
for $0<h\le h_0(\tau,A)\ll 1$, with a constant $C>0$ independent of $h$, $\tau$ and $A$. 
\end{lemma}

We are going to use (3.1) with $f=\left(P_\varphi^L(h)-i\right)^{-1}g$. In view of the identity
$$1=(1-i\mp i\varepsilon)\left(P_\varphi^L(h)-i\right)^{-1}+\left(P_\varphi^L(h)-i\right)^{-1}
\left(P_\varphi^L(h)-1\pm i\varepsilon\right)$$
and Lemma 3.2, we have
$$\left\|\langle x\rangle^{-s}g\right\|_{H^1}\le 2\left\|\langle x\rangle^{-s}
\left(P_\varphi^L(h)-i\right)^{-1}g\right\|_{H^1}$$
$$+\left\|\langle x\rangle^{-s}\left(P_\varphi^L(h)-i\right)^{-1}
\left(P_\varphi^L(h)-1\pm i\varepsilon\right)g\right\|_{H^1}$$
 $$\le\frac{2C_1}{h}\left\|\mu^{-1}\left(P_\varphi^L(h)-i\right)^{-1}
\left(P_\varphi^L(h)-1\pm i\varepsilon\right)g\right\|_{L^2}
+2C_2\left(\frac{\varepsilon}{h}\right)^{1/2}\left\|\left(P_\varphi^L(h)-i\right)^{-1}g\right\|_{H^1}$$
$$+C_3\left\|\left(P_\varphi^L(h)-i\right)^{-1}\left(P_\varphi^L(h)-1\pm i\varepsilon\right)g\right\|_{H^1}$$
$$\le\frac{2C_1}{h}\left\|\mu^{-1}\left(P_\varphi^L(h)-i\right)^{-1}\mu\right\|_{H^{-1}\to L^2}\left\|\mu^{-1}
\left(P_\varphi^L(h)-1\pm i\varepsilon\right)g\right\|_{H^{-1}}$$
 $$+2C_2\left(\frac{\varepsilon}{h}\right)^{1/2}\left\|\left(P_\varphi^L(h)-i\right)^{-1}\right\|_{L^2\to H^1}\|g\|_{L^2}
$$ $$+C_3\left\|\left(P_\varphi^L(h)-i\right)^{-1}\right\|_{H^{-1}\to H^1}\left\|\left(P_\varphi^L(h)-
1\pm i\varepsilon\right)g\right\|_{H^{-1}}$$
$$\le\frac{C'_1}{h}\left\|\mu^{-1}\left(P_\varphi^L(h)-1\pm i\varepsilon\right)g\right\|_{H^{-1}}
+C'_2\left(\frac{\varepsilon}{h}\right)^{1/2}\|g\|_{L^2}+C'_3\left\|\left(P_\varphi^L(h)-
1\pm i\varepsilon\right)g\right\|_{H^{-1}}\eqno{(3.15)}$$
with a constant $C'_1>0$ independent of $\varepsilon$, $h$, $\tau$, $A$ and $g$, 
and constants $C'_2, C'_3>0$ independent of $\varepsilon$, $h$ and $g$. Since the function $\mu$ is bounded on ${\bf R}^n$,
there exists $0<h_0(\varphi)\ll 1$ such that for $0<h\le h_0$ the last term in the right-hand side of (3.15)
can be bounded by the first one. Thus we get (2.1) from (3.15).

\section{Proof of Lemma 3.2}

The estimate (3.14) is known to hold with $\ell=0$ and $P^L_\varphi(h)$ replaced by $-h^2\Delta$. We will use this
to show that (3.14) with $\ell=0$ still holds for first-order perturbations of the form $-h^2\Delta+Q(h)$, where
$$Q(h)=\sum_{|\alpha|=1}q_\alpha^{(1)}(x,h){\cal D}_x^\alpha+\sum_{|\alpha|=1}{\cal D}_x^\alpha q_\alpha^{(2)}(x,h)
+q_0(x,h)$$
with coefficients satisfying
$$\left|q_\alpha^{(1)}(x,h)\right|+\left|q_\alpha^{(2)}(x,h)\right|+\left|q_0(x,h)\right|\le Ch,\quad\forall x\in{\bf R}^n.
\eqno{(4.1)}$$
Clearly, (4.1) implies
$$\|Q(h)\|_{H^1\to H^{-1}}\le Ch.\eqno{(4.2)}$$
By (4.2) and the resolvent identity
$$\left(-h^2\Delta+Q(h)-i\right)^{-1}=\left(-h^2\Delta-i\right)^{-1}+\left(-h^2\Delta-i\right)^{-1}Q(h)
\left(-h^2\Delta+Q(h)-i\right)^{-1}$$
we get
$$\left\|\left(-h^2\Delta+Q(h)-i\right)^{-1}\right\|_{H^{-1}\to H^1}\le\left\|\left(-h^2\Delta-
i\right)^{-1}\right\|_{H^{-1}\to H^1}$$ $$+\left\|\left(-h^2\Delta-i\right)^{-1}\right\|_{H^{-1}\to H^1}
\|Q(h)\|_{H^1\to H^{-1}}\left\|\left(-h^2\Delta+Q(h)-i\right)^{-1}\right\|_{H^{-1}\to H^1}$$
$$\le C+O(h)\left\|\left(-h^2\Delta+Q(h)-i\right)^{-1}\right\|_{H^{-1}\to H^1}.\eqno{(4.3)}$$
Now, taking $h$ small enough (depending on the coefficients of $Q(h)$) we can absorb the last term in the 
right-hand side of (4.3) and obtain the desired estimate with a constant $C>0$ independent of $q_\alpha^{(1)}$,
$q_\alpha^{(2)}$, $q_0$ and $h$.

Thus, to prove (3.14) it suffices to show that the operator $\mu^{-\ell}P_\varphi^L(h)\mu^\ell$ is equal to
$-h^2\Delta$ plus a first-order differential operator with coefficients satisfying (4.1). To do so, observe 
first that $\mu^{-\ell}P_\varphi^L(h)\mu^\ell=P_\psi^L(h)$, where $\psi=\varphi-\ell\log\mu$. Furthermore, we have
$$P_\psi^L(h)=-h^2\Delta+\left(i\widetilde b^L-h\nabla\psi\right)\cdot h\nabla+h\nabla\cdot
\left(i\widetilde b^L-h\nabla\psi\right)-h^2|\nabla\psi|^2-2ih\widetilde b^L\cdot\nabla\psi+\widetilde V^L.$$
It is easy to see that $|\psi'(r)|$ is bounded on ${\bf R}$, and hence 
$|\nabla\psi(|x|)|$ is bounded on ${\bf R}^n$. This together with the assumptions on $\widetilde b^L$
and $\widetilde V^L$ imply the desired properties of the coefficients of the operator $P_\psi^L(h)$. 
\eproof

G. Vodev

Universit\'e de Nantes,

 D\'epartement de Math\'ematiques, UMR 6629 du CNRS,
 
 2, rue de la Houssini\`ere, BP 92208, 
 
 44332 Nantes Cedex 03, France,
 
 e-mail: vodev@math.univ-nantes.fr

\end{document}